\documentclass[a4paper,12pt]{amsart}
\usepackage{amssymb}
\usepackage{ifthen}
\usepackage{graphicx}
\nonstopmode \numberwithin{equation}{section}
\usepackage[usenames]{color}
\newtheorem{thm}{Theorem}[section]
\newtheorem{cor}{Corollary}[section]%[equation]
\newtheorem{lem}{Lemma}[section]

\newtheorem{rem}{Remark}[section]%[equation]
\newtheorem{rems}[equation]{Remarks}
\theoremstyle{definition}
\newtheorem{defin}{Definition}[section]%[equation]
\newtheorem{examp}[equation]{Example}%[section]
\newtheorem{prob}[equation]{Problem}
\newtheorem{ques}[equation]{Question}
\newtheorem{op}{Open Problem}[section]%[equation]

\newtheorem{conj}[equation]{Conjecture}
\newtheorem{deter}[equation]{Determination}
\newtheorem{case}{Case}[section]%[equation]
\newtheorem{subcase}[equation]{Subcase}
\newtheorem{claim}{Claim}[section]%[equation]
\newtheorem{subclaim}{Subclaim}

\newcounter {own}
\def\theown {\thesection       .\arabic{own}}

\newenvironment{pf}[1][]{%
 \vskip 3mm
 \noindent
 \ifthenelse{\equal{#1}{}}%
  {{\bf Proof. }}%
  {{\bf #1.} }%
 }%
{\qed\bigskip}

\newcounter{alphabet}
\newcounter{tmp}
\newenvironment{Thm}[1][]{\refstepcounter{alphabet}%
\bigskip%
\noindent%
{\bf Theorem \Alph{alphabet}}%
\ifthenelse{\equal{#1}{}}{}{ (#1)}%
{\bf .} \itshape}{\vskip 8pt}
\makeatletter
\newcommand{\Ref}[1]{\@ifundefined{r@#1}{}{\setcounter{tmp}{\ref{#1}}\Alph{tmp}}}
\makeatother

\newcommand{\IR}{{\mathbb R}}

\newcommand{\diam}{{\operatorname{diam}}}

\def\be{\begin{equation}}
\def\ee{\end{equation}}

\newcommand{\bee}{\begin{enumerate}}
\newcommand{\eee}{\end{enumerate}}

\newcommand{\blem}{\begin{lem}}
\newcommand{\elem}{\end{lem}}
\newcommand{\bthm}{\begin{thm}}
\newcommand{\ethm}{\end{thm}}
\newcommand{\bcor}{\begin{cor}}
\newcommand{\ecor}{\end{cor}}
\newcommand{\beg}{\begin{examp}}
\newcommand{\eeg}{\end{examp}}
\newcommand{\begs}{\begin{examples}}
\newcommand{\eegs}{\end{examples}}
\newcommand{\bdefe}{\begin{defin}}
\newcommand{\edefe}{\end{defin}}
\newcommand{\bprob}{\begin{prob}}
\newcommand{\eprob}{\end{prob}}
\newcommand{\bques}{\begin{ques}}
\newcommand{\eques}{\end{ques}}
\newcommand{\bei}{\begin{itemize}}
\newcommand{\eei}{\end{itemize}}

\newcommand{\bde}{\begin{deter}}
\newcommand{\ede}{\end{deter}}
\newcommand{\bca}{\begin{case}}
\newcommand{\eca}{\end{case}}
\newcommand{\bsca}{\begin{subcase}}
\newcommand{\esca}{\end{subcase}}
\newcommand{\bcl}{\begin{claim}}
\newcommand{\ecl}{\end{claim}}
\newcommand{\bscl}{\begin{subclaim}}
\newcommand{\escl}{\end{subclaim}}
\newcommand{\bcon}{\begin{conj}}
\newcommand{\econ}{\end{conj}}
\newcommand{\bcons}{\begin{conjs}}
\newcommand{\econs}{\end{conjs}}
\newcommand{\bprop}{\begin{propo}}
\newcommand{\eprop}{\end{propo}}
\newcommand{\br}{\begin{rem}}
\newcommand{\er}{\end{rem}}
\newcommand{\brs}{\begin{rems}}
\newcommand{\ers}{\end{rems}}
\newcommand{\bo}{\begin{obser}}
\newcommand{\eo}{\end{obser}}
\newcommand{\bos}{\begin{obsers}}
\newcommand{\eos}{\end{obsers}}
\newcommand{\bpf}{\begin{pf}}
\newcommand{\epf}{\end{pf}}
\newcommand{\ba}{\begin{array}}
\newcommand{\ea}{\end{array}}
\newcommand{\beq}{\begin{eqnarray}}
\newcommand{\beqq}{\begin{eqnarray*}}
\newcommand{\eeq}{\end{eqnarray}}
\newcommand{\eeqq}{\end{eqnarray*}}

\newcommand{\ds}{\displaystyle}

\newcommand{\bop}{\begin{op}}
\newcommand{\eop}{\end{op}}

\newtheorem{pfofThm1.5}[equation]{}

\newcounter{minutes}
\divide\time by 60
\newcounter{hours}
\multiply\time by 60 \addtocounter{minutes}{-\time}
\begin{document}

\bibliographystyle{amsplain}

\title{On quasim\"obius maps in real Banach spaces}

%%%%%%%% BEGIN TIMESTAMP
\def\thefootnote{}
\footnotetext{ \texttt{\tiny File:~\jobname .tex,
          printed: \number\year-\number\month-\number\day,
          \thehours.\ifnum\theminutes<10{0}\fi\theminutes}
} \makeatletter\def\thefootnote{\@arabic\c@footnote}\makeatother
%%%%%%%% END TIMESTAMP

\author{M. Huang}
\address{M. Huang, Department of Mathematics,
Hunan Normal University, Changsha,  Hunan 410081, People's Republic
of China} \email{mzhuang79@yahoo.com.cn}

\author{Y. Li}
\address{Y. Li, Department of Mathematics,
Hunan Normal University, Changsha,  Hunan 410081, People's Republic
of China} \email{yaxiangli@163.com}

\author{M. Vuorinen}
\address{M. Vuorinen, Department of Mathematics and Statistics,
University of Turku, 20014 Turku, Finland } \email{vuorinen@utu.fi}

\author{X. Wang $^* $
%${}^{~\mathbf{*}}$
}
\address{X. Wang, Department of Mathematics,
Hunan Normal University, Changsha,  Hunan 410081, People's Republic
of China} \email{xtwang@hunnu.edu.cn}

\date{}
 \keywords{
Uniform domain, quasim\"obius map, quadruple, CQH homeomorphism,
neargeodesic, solid arc.\\
${}^{\mathbf{*}}$ Corresponding author}

\begin{abstract}
Suppose that $E$ and $E'$ denote real Banach spaces with dimension
at least $2$, that $D\varsubsetneq E$ and $D'\varsubsetneq E'$ are
domains, that $f: D\to D'$ is an $(M,C)$-CQH homeomorphism, and that
$D$ is uniform.
 The main aim of this paper is to prove that $D'$ is a uniform
domain if and only if $f$ extends to a homeomorphism $\overline{f}:
\overline{D}\to \overline{D}'$ and $\overline{f}$ is $\eta$-QM
relative to $\partial D$. This result shows that the answer to one
of the open problems raised by V\"ais\"al\"a from 1991 is
affirmative.
\end{abstract}

\thanks{The research was partly supported by NSFs of
China (No. 11071063 and No. 11101138) and by the Academy of Finland,
Project 2600066611.}

\maketitle{} \pagestyle{myheadings} \markboth{}{On quasim\"obius
maps in real Banach spaces}

\section{Introduction and main results}\label{sec-1}

During the past three decades, the quasihyperbolic metric has become
an important tool in geometric function theory and in its
generalizations to metric spaces and to Banach spaces \cite{Vai5}. Yet,
some basic questions of the quasihyperbolic geometry in Banach
spaces are open. For instance, only recently the convexity of
quasihyperbolic balls has been studied in \cite{k,krt,rt,Vai9} in the setup
of Banach spaces.

Our study is motivated by V\"ais\"al\"a's theory of freely
quasiconformal mappings  and other related maps in the setup of
Banach spaces \cite{Vai6-0, Vai6, Vai5}. Our goal is to study some
of the open problems formulated by him. We begin with some basic
definitions and the statements of our results. The proofs and
necessary supplementary notation and terminology will be given
thereafter.

Throughout the paper, we always assume that $E$ and $E'$ denote real
Banach spaces with dimension at least $2$. The norm of a vector $z$
in $E$ is written as $|z|$, and for every pair of points $z_1$,
$z_2$ in $E$, the distance between them is denoted by $|z_1-z_2|$,
the closed line segment with endpoints $z_1$ and $z_2$ by $[z_1,
z_2]$. We begin with the following concepts following closely the
notation and terminology of \cite{TV, Vai2, Vai, Vai6-0, Vai6} or
\cite{Martio-80}.

\bdefe \label{def1.3} A domain $D$ in $E$ is called $c$-{\it
uniform} in the norm metric provided there exists a constant $c$
with the property that each pair of points $z_{1},z_{2}$ in $D$ can
be joined by a rectifiable arc $\alpha$ in $ D$ satisfying

 \bee
\item\label{wx-4} $\ds\min_{j=1,2}\ell (\alpha [z_j, z])\leq c\,d_{D}(z)
$ for all $z\in \alpha$, and

\item\label{wx-5} $\ell(\alpha)\leq c\,|z_{1}-z_{2}|$,
\eee

\noindent where $\ell(\alpha)$ denotes the length of $\alpha$,
$\alpha[z_{j},z]$ the part of $\alpha$ between $z_{j}$ and $z$, and
$d_D(z)$ the distance from $z$ to the boundary $\partial D$ of $D$.
Also, we say that $\gamma$ is a double $c$-cone arc. \edefe

In \cite{Vai2}, V\"ais\"al\"a obtained the following result
concerning the relation between the class of uniform domains and
quasim\"obius (briefly, QM) maps (see Definition \ref{def2'}) in $\overline{\mathbb R}^n = {\mathbb R}^n \cup
\{\infty\}$.

\begin{Thm}\label{ThmA-1}$($\cite[Theorem 5.6]{Vai2}$)$
Suppose that $n\geq 2$, that $D$ is a
 $c$-uniform domain in $\overline{\mathbb{R}}^n$ and that $f:D\to D'$
is a quasiconformal mapping. Then the following conditions are
quantitatively equivalent:

$(1)$ $D'$ is a $c_1$-uniform domain;

$(2)$ $f$ is $\eta$-QM.
\end{Thm}

In \cite{Vai6}, V\"ais\"al\"a generalized Theorem \Ref{ThmA-1} to
the case of Banach spaces. His result is as follows.

\begin{Thm}\label{ThmA}$($\cite[Theorem 7.18]{Vai6}$)$
Let $D$ and $D'$ be domains in $E$ and $E'$, respectively. Suppose
that $D$ is a $c$-uniform domain and that $f:D\to D'$ is
$\varphi$-FQC $($see Definition \ref{def1.7-2}$)$. Then the
following conditions are quantitatively equivalent:

$(1)$ $D'$ is a $c_1$-uniform domain;

$(2)$ $f$ is $\eta$-QM.
\end{Thm}

Further, V\"ais\"al\"a \cite[7.19]{Vai6} raised the following open
problem.

\bop\label{Con1} Does Theorem \Ref{ThmA} remain true if
$\varphi$-FQC is replaced by $(M,C)$-CQH $($see
Definition \ref{def1.6}$)$ and $\eta$-QM by $\eta$-QM rel $\partial D$,
respectively? \eop

Studying this problem, V\"ais\"al\"a proved the following result.

\begin{Thm}\label{ThmC}$($\cite[Theorem 7.9]{Vai6}$)$
Suppose that $D\varsubsetneq E$ and $D'\varsubsetneq E'$ are
$c$-uniform domains and that $f:D\to D'$ is $(M,C)$-CQH. Then $f$
extends to a homeomorphism $\overline{f}: \overline{D}\to
\overline{D'}$ and $\overline{f}$ is $\eta$-QM rel $\partial D$ with
$\eta$ depending only on $(M,C,c)$. In particular,
$\overline{f}|\partial D$ is $\eta$-QM.\end{Thm}

The aim of this paper is to discuss Open Problem \ref{Con1} further.
Our main result is the next theorem, which shows that the answer to Open
Problem \ref{Con1} is in the affirmative.

\begin{thm}\label{thm1.1}
Suppose that $D$ is a $c$-uniform domain and that $f: D\to D'$ is
$(M,C)$-CQH, where $D\varsubsetneq E$ and $D'\varsubsetneq E'$. Then
the following conditions are quantitatively equivalent:

$(1)$ $D'$ is a $c_1$-uniform domain;

$(2)$  $f$ extends to a homeomorphism $\overline{f}: \overline{D}\to
\overline{D}'$ and $\overline{f}$ is $\eta$-QM rel $\partial D$.
\end{thm}

The organization of this paper is as follows. Section \ref{sec-2} contains some preliminaries and a new lemma. The proof of the main result is given in Section \ref{sec-3}.

%%%%%%%%%%%%%%%%%
%%%%%%%%%%%%%%%%%
\section{Preliminaries}\label{sec-2}
%%%%%%%%%%%%%%%%%
%%%%%%%%%%%%%%%%%

\subsection{Quasihyperbolic distance, quasihyperbolic geodesics, neargeodesics and solid arcs}

The {\it quasihyperbolic length} of a rectifiable arc or a path
$\alpha$ in the norm metric in $D$ is the number (cf.
\cite{GP,Vai3}):

$$\ell_k(\alpha)=\int_{\alpha}\frac{|dz|}{d_{D}(z)}.
$$

For each pair of points $z_1$, $z_2$ in $D$, the {\it quasihyperbolic distance}
$k_D(z_1,z_2)$ between $z_1$ and $z_2$ is defined in the usual way:
$$k_D(z_1,z_2)=\inf\ell_k(\alpha),
$$
where the infimum is taken over all rectifiable arcs $\alpha$
joining $z_1$ to $z_2$ in $D$. For all $z_1$, $z_2$ in $D$, we have
(cf. \cite{Vai3})

\beq\label{eq(0000)} k_{D}(z_1, z_2)\geq
\inf\left\{\log\Big(1+\frac{\ell(\alpha)}{\min\{d_{D}(z_1),
d_{D}(z_2)\}}\Big)\right\}\geq \Big|\log
\frac{d_{D}(z_2)}{d_{D}(z_1)}\Big|,\eeq where the infimum is taken
over all rectifiable curves $\alpha$ in $D$ connecting $z_1$ and
$z_2$.  Moreover, if $|z_1-z_2|\le d_D(z_1)$, we have
\cite{Vai6-0, Mvo1}
\begin{equation} \label{upperbdk}
k_D(z_1,z_2)\le \log\Big( 1+ \frac{
|z_1-z_2|}{d_D(z_1)-|z_1-z_2|}\Big).
\end{equation}

The quasihyperbolic metric of a domain in $\IR^n$ was introduced by Gehring and Palka \cite{GP} and it has been recently used by many authors
 in the study of quasiconformal mappings and related questions; see
\cite{Avv, Bea, BHX, Geo, HIMPS, KN, Vai6-0, Vai6, Vai5} etc.

In \cite{Vai6},  V\"ais\"al\"a characterized uniform domains by the quasihyperbolic metric.

\begin{Thm}\label{thm0.1} {\rm (\cite[Theorem 6.16]{Vai6})}
For a domain $D$ in $E$, the following are quantitatively
equivalent: \bee

\item $D$ is a $c$-uniform domain;
\item $k_D(z_1,z_2)\leq c'\;
 \log\left(1+\frac{\ds{|z_1-z_2|}}{\ds\min\{d_{D}(z_1),d_{D}(z_2)\}}\right)$ for all $z_1,z_2\in D$;
\item $k_D(z_1,z_2)\leq c'_1\;
 \log\left(1+\frac{\ds{|z_1-z_2|}}{\ds\min\{d_{D}(z_1),d_{D}(z_2)\}}\right)+d$ for all $z_1,z_2\in D$.\eee
\end{Thm}
\begin{rem}\label{thm0-1'} By \cite[Theorem 2.23]{Vai4}, we know that in Theorem \Ref{thm0.1}, if $(1)$ holds,
then $(2)$ holds with $c'\leq 7c^3.$
\end{rem}

 In the case of domains in $ {\mathbb R}^n \,,$ the equivalence
  of items (1) and (3) in Theorem D is due to Gehring and Osgood \cite{Geo} and the
  equivalence of items (2) and (3) due to Vuorinen \cite{Mvo1}.

Recall that an arc $\alpha$ from $z_1$ to
$z_2$ is a {\it quasihyperbolic geodesic} if
$\ell_k(\alpha)=k_D(z_1,z_2)$. Each subarc of a quasihyperbolic
geodesic is obviously a quasihyperbolic geodesic. It is known that a
quasihyperbolic geodesic between every pair of points in $E$ exists if the
dimension of $E$ is finite, see \cite[Lemma 1]{Geo}. This is not
true in arbitrary spaces (cf. \cite[Example 2.9]{Vai6-0}).
In order to remedy this shortage, V\"ais\"al\"a introduced the following concepts \cite{Vai6}.

\bdefe \label{def1.4}
 Let $\alpha$ be an arc in $E$. The arc
may be closed, open or half open. Let $\overline{x}=(x_0,...,x_n)$,
$n\geq 1$, be a finite sequence of successive points of $\alpha$.
For $h\geq 0$, we say that $\overline{x}$ is {\it $h$-coarse} if
$k_D(x_{j-1}, x_j)\geq h$ for all $1\leq j\leq n$. Let $\Phi_{k}(\alpha,h)$
be the family of all $h$-coarse sequences of $\alpha$. Set

$$s_{k}(\overline{x})=\sum^{n}_{j=1}k_D(x_{j-1}, x_j)$$ and
$$\ell_{k}(\alpha, h)=\sup \{s_{k}(\overline{x}): \overline{x}\in \Phi_{k}(\alpha,h)\}$$
with the agreement that $\ell_{k}(\alpha, h)=0$ if
$\Phi_{k}(\alpha,h)=\emptyset$. Then the number $\ell_{k}(\alpha, h)$ is the
{\it $h$-coarse quasihyperbolic length} of $\alpha$.\edefe

\bdefe \label{def1.5} Let $D$ be a domain in $E$. An arc $\alpha\subset D$
is {\it $(\nu, h)$-solid} with $\nu\geq 1$ and $h\geq 0$ if
$$\ell_k(\alpha[x,y], h)\leq \nu\;k_D(x,y)$$ for all $x, y\in \alpha$.
A {\it $(\nu,0)$-solid arc} is said to be a {\it $\nu$-neargeodesic}, i.e.
an arc $\alpha\subset D$ is a $\nu$-neargeodesic if and only if $\ell_k(\alpha[x,y])\leq \nu\;k_D(x,y)$
for all $x, y\in \alpha$.\edefe

Obviously, a $\nu$-neargeodesic is a quasihyperbolic geodesic if and
only if $\nu=1$.

In \cite{Vai6}, V\"ais\"al\"a established the following property concerning
the existence of neargeodesics in $E$.

\begin{Thm}\label{LemA} $($\cite[Theorem 3.3]{Vai6}$)$
Let $\{z_1,\, z_2\}\subset D$ and $\nu>1$. Then there is a
$\nu$-neargeodesic in $D$ joining $z_1$ and $z_2$.
\end{Thm}

The following result due to V\"ais\"al\"a is from \cite{Vai6}.

\begin{Thm}\label{Thm4-1} {\rm (\cite[Theorem 6.22]{Vai6})} Suppose that
$\gamma\subset G\varsubsetneq E$ is a $(\nu,h)$-solid arc with
endpoints $a_0$, $a_1$ and that $G$ is a $c$-uniform domain. Then
there is a constant $c_2=c_2(\nu,h,c)\geq 1$ such that

\begin{enumerate}
\item  $\ds\min\Big\{\diam(\gamma [a_0, z]),\diam(\gamma [a_1, z])\Big\}\leq c_2\,d_{G}(z)
$ for all $z\in \gamma$, and

\item  $\diam(\gamma)\leq c_2\max\Big\{|a_0-a_1|,2(e^h-1)\min\{d_{G}(a_0),d_{G}(a_1)\}\Big\}$.\end{enumerate}
 \end{Thm}

\subsection{Quasisymmetric homeomorphisms and quasim\"obius maps}

Let $X$ be a metric space and $\dot{X}=X\cup \{\infty\}$. By a
triple in $X$ we mean an ordered sequence $T=(x,a,b)$ of three
distinct points in $X$. The ratio of $T$ is the number
$$\rho(T)=\frac{|a-x|}{|b-x|}.$$ If $f: X\to Y$ is  an injective
map, the image of a triple  $T=(x,a,b)$  is the triple
$fT=(fx,fa,fb)$.

Suppose that  $A\subset X$. A triple  $T=(x,a,b)$ in $X$ is said to
be a triple in the pair $(X, A)$ if $x\in A$ or if $\{a,b\}\subset
A$. Equivalently, both $|a-x|$ and $|b-x|$ are distances from a
point in $A$.

\bdefe \label{def1-0} Let $X$ and $Y$ be two metric spaces, and let
$\eta: [0, \infty)\to [0, \infty)$ be a homeomorphism. Suppose
$A\subset X$. An embedding $f: X\to Y$ is said to be {\it
$\eta$-quasisymmetric} relative to $A$, or briefly $\eta$-$QS$ rel
$A$, if $\rho(f(T))\leq \eta(\rho(T))$  for each triple $T$ in
$(X,A)$. \edefe

It is known that an embedding $f: X\to Y$ is $\eta$-$QS$ rel $A$ if
and only if $\rho(T)\leq t$ implies that $\rho(f(T))\leq \eta(t)$
for each triple $T$ in $(X,A)$ and $t\geq 0$ (cf. \cite{TV}).
Obviously, ``quasisymmetric rel $X$" is equivalent to ordinary
``quasisymmetric".

A quadruple in $X$ is an ordered sequence $Q=(a,b,c,d)$ of four
distinct points in $X$. The cross ratio of $Q$ is defined to be the
number
$$\tau(Q)=|a,b,c,d|=\frac{|a-b|}{|a-c|}\cdot\frac{|c-d|}{|b-d|}.$$ Observe that the definition is extended in
the well known manner to the case where one of the points is
$\infty$. For example,
$$|a,b,c,\infty|= \frac{|a-b|}{|a-c|}.$$ If $X_0 \subset \dot{X}$ and if $f: X_0\to \dot{Y}$
is an injective map, the image of a quadruple $Q$ in $X_0$ is the
quadruple $fQ=(fa,fb,fc,fd)$. Suppose that $A\subset X_0$. We say
that a quadruple $Q=(a,b,c,d)$ in $X_0$ is a quadruple in the pair
$(X_0, A)$ if $\{a,d\}\subset A$ or $\{b,c\}\subset A$.
Equivalently, all four distances in the definition of $\tau(Q)$ are
(at least formally) distances from a point in $A$.

\bdefe \label{def2'} Let $\dot{X}$ and $\dot{Y}$ be two metric
spaces and let $\eta: [0, \infty)\to [0, \infty)$ be a
homeomorphism. Suppose $A\subset \dot{X}$. An embedding $f:
\dot{X}\to \dot{Y}$ is said to be {\it $\eta$-quasim\"obius}
relative to $A$, or briefly $\eta$-$QM$ rel $A$, if the inequality
$\tau(f(Q))\leq \eta(\tau(Q))$ holds for each quadruple in $(X,A)$.
\edefe

Apparently, ``$\eta$-$QM$ rel $X$" is equivalent to ordinary
``quasim\"obius".

\subsection{Coarsely quasihyperbolic homeomorphisms and freely
quasiconformal mappings}

\bdefe \label{def1.6} We say that a homeomorphism $f: D\rightarrow
D'$ is {\it $C$-coarsely $M$-quasihyperbolic}, or briefly
$(M,C)$-CQH, in the quasihyperbolic metric if it satisfies
$$\frac{k_D(x,y)-C}{M}\leq k_{D'}(f(x),f(y))\leq M\;k_D(x,y)+C$$
for all $x$, $y\in D$. \edefe

The following result shows that the class of solid arcs is invariant under the CQH homeomorphisms.

\begin{Thm}\label{LemB2} {\rm (\cite[Theorem 4.15]{Vai6})} For domains
$D\varsubsetneq E$ and $D'\varsubsetneq E'$, suppose that $f: D\to
D'$ is $(M,C)$-CQH. If $\gamma$ is a $(\nu_1,h_1)$-solid arc in $D$,
then the arc $f(\gamma)$ is $(\nu,h)$-solid in $D'$ with $(\nu, h)$
depending only on $(\nu_1, h_1, M, C)$.\end{Thm}

\bdefe \label{def1.7-2} Let $G\not=E$ and $G'\not=E'$ be metric
spaces, and let $\varphi:[0,\infty)\to [0,\infty)$ be a growth
function, that is, a homeomorphism with $\varphi(t)\geq t$. We say
that a homeomorphism $f: G\to G'$ is {\it $\varphi$-semisolid} if
$$ k_{G'}(f(x),f(y))\leq \varphi(k_{G}(x,y))$$
for all $x$, $y\in G$, and {\it $\varphi$-solid} if both $f$ and
$f^{-1}$ satisfy this condition.

We say that $f$ is {\it fully $\varphi$-semisolid} (resp. {\it fully
$\varphi$-solid}) if $f$ is $\varphi$-semisolid (resp.
$\varphi$-solid) on every  subdomain of $G$. In particular, when
$G=E$, the corresponding subdomains are taken to be proper ones.
Fully $\varphi$-solid mappings are also called {\it freely
$\varphi$-quasiconformal mappings}, or briefly {\it $\varphi$-FQC
mappings}.\edefe

\subsection{Basic assumptions and a lemma}
For convenience, in the following, we always assume that $x$, $y$, $z$, $\ldots$
denote points in $D$ and $x'$, $y'$, $z'$, $\ldots$
the images in $D'$ of $x$, $y$, $z$, $\ldots$
under $f$, respectively. Also we assume that $\alpha$, $\beta$, $\gamma$, $\ldots$
denote curves in $D$ and $\alpha'$, $\beta'$, $\gamma'$, $\ldots$  the images in $D'$ of
$\alpha$, $\beta$, $\gamma$, $\ldots$
under $f$, respectively.
\medskip

\noindent {\bf Basic assumption A. }\quad Let $G$ be a domain in
$E$. For $x$, $y\in G$, let $\beta$ be a $2$-neargeodesic joining
$x$ and $y$ in $G$. Suppose that $G'$ is a $c$-uniform domain in
$E'$ and $f: G\to G'$ is an $(M,C)$-CQH homeomorphism. It follows
from Theorem \Ref{LemB2} that $\beta'$ is $(\nu,h)$-solid, where
$(\nu,h)$ depends only on $(M, C, c)$. Without loss of generality,
we may assume that $d_{G'}(y')\geq d_{G'}(x')$. Then there must
exist a point $z'_0\in\beta'$ which is the first point on $\beta'$ in the
direction from $x'$ to $y'$ such that

$$d_{G'}(z'_0)=\sup\limits_{p'\in \beta'}d_{G'}(p').
$$
It is possible that $z'_0=x'$ or $y'$. Then we have

\begin{lem} \label{lem-j-j}
\noindent $(1)$ For all  $z'\in
\beta'[x', z'_0]$,
 $$|x'-z'|\leq
\mu_1\;d_{G'}(z'),$$ and for all  $z'\in \beta'[y', z'_0]$,
$$|y'-z'|\leq \mu_1\;d_{G'}(z');$$

\noindent $(2)$ $\diam(\beta')\leq
\mu_1\max\{|x'-y'|,2(e^{h}-1)d_{G'}(x')\},$

\noindent where $\mu_1=4c_2^2$, $c_2=c_2(\nu,h,c)$ is the same as in Theorem \Ref{Thm4-1},
and $\nu$, $h$ and $c$ are as in Basic assumption A.
\end{lem}

\bpf By Theorem \Ref{Thm4-1}, it suffices to prove the first assertion in $(1)$.
For the case $\min\{\diam(\beta'[x',z']),\diam(\beta'[y',z'])\}=\diam(\beta'[x',z'])$,
it follows from Theorem \Ref{Thm4-1} that the proof is obvious. For the other case
$\min\{\diam(\beta'[x',z']),\diam(\beta'[y',z'])\}=\diam(\beta'[y',z'])$, we first have the following claim.

\bcl\label{xt-1} $\diam(\beta'[x',z'])\leq 2c_2d_{G'}(z'_0)$.\ecl

Suppose on the contrary that $$\diam(\beta'[x',z'])>
2c_2d_{G'}(z'_0).$$
Obviously, there must exist some point $w'\in\beta'[x',z']$ such that
$$\diam(\beta'[w',z'])=\frac{1}{2}\diam(\beta'[x',z'])\;\, \mbox{and}\;\, \diam(\beta'[x',w'])\geq\frac{1}{2}\diam(\beta'[x',z']).$$
It follows from Theorem \Ref{Thm4-1} that
\begin{eqnarray*}c_2d_{G'}(w')&\geq&\min\{\diam(\beta'[x',w']),\diam(\beta'[y',w'])\}
\\ \nonumber&\geq& \frac{1}{2}\,\diam(\beta'[x',z'])>c_2d_{G'}(z'_0 ). \end{eqnarray*}
This is the desired contradiction which completes the proof of Claim
\ref{xt-1}.
\medskip

If $\diam(\beta'[y',z'])\leq \frac{1}{2}d_{G'}(z'_0)$, then by Claim
\ref{xt-1},
$$|x'-z'|\leq\diam(\beta'[x',z'])
\leq 2c_2\,d_{G'}(z'_0)\leq
4c_2d_{G'}(z'),$$ since $d_{G'}(z')\geq
d_{G'}(z'_0)-|z'_0-z'|$.

If $\diam(\beta'[y',z'])> \frac{1}{2}d_{G'}(z'_0)$, then we see from Claim
\ref{xt-1} and
Theorem
\Ref{Thm4-1} that
$$|x'-z'|\leq\diam(\beta'[x',z'])
\leq 2c_2\,d_{G'}(z'_0)\leq
4c_2\diam(\beta'[y',z'])\leq4c_2^2
d_{G'}(z').$$
The proof is finished.\epf

%%%%%%%%%%%%%%%%%
%%%%%%%%%%%%%%%%%
\section{The proof of Theorem \ref{thm1.1}} \label{sec-3}
%%%%%%%%%%%%%%%%%%%%%%%%%%%%%%%%%%%%%%
%%%%%%%%%%%%%%%%%%%%%%%%%%%%%%%%%%%%%%
First, we recall the following results which are from  \cite{Vai2} and \cite{Vai6}, respectively.

\begin{Thm}\label{LemC} $($\cite[Theorem 3.19]{Vai2}$)$
Suppose that $A\subset\dot{X}$,
that $f: A \to \dot{X}$ is $\eta$-QM,
and suppose that $\overline{f(A)}\backslash \{\infty\}$
is complete. Then $f$ has a unique extension to an $\eta$-QM embedding $g: \overline{A}\to \dot{X}$.
\end{Thm}

\begin{Thm}\label{LemB} $($\cite[Theorem 6.26]{Vai6}$)$
Suppose that $f:D \to D'$ is $\eta$-QM
and that $D$ is a $c$-uniform domain. Then $D'$ is a $c_1$-uniform
domain, where $c_1$ depends only on $c$ and $\eta$.
\end{Thm}

For more details of quasisymmetric homeomorphisms and quasi\-m\"obius maps, the reader is referred to \cite{TV, Vai2, Vai6-0}.

The proof of Theorem \ref{thm1.1} will be accomplished through a
series of lemmas. Before the statements of the lemmas, we give
another basic assumption.
\medskip

\noindent {\bf Basic assumption B.}\quad
Throughout this
section, we always assume that $D$ is a $c$-uniform domain, that $f:
D\to D'$ is $(M,C)$-CQH, where $D\varsubsetneq E$ and
$D'\varsubsetneq E'$, and that $f$ extends to a homeomorphism
$\overline{f}: \overline{D}\to \overline{D}'$ and $\overline{f}$ is
$\eta$-QM rel $\partial D$. By auxiliary translations and
inversions, it follows from Theorems \Ref{LemB} and \Ref{LemC} that
we may normalize the map $f$ and the domain $D$ so that $\infty\in
\partial D$ and $\overline{f}(\infty)=\infty$.
 Then $f$ is $\eta$-QS rel
$\partial D$.
%\medskip

\noindent {\bf Constants.}\quad
For the convenience of the statements of the lemmas below, we write down the related constants:

\noindent $(1)$ $\mu_2=\max\{4(e^{h}-1)\mu_1, 6\mu_1\}$,

\noindent $(2)$ $\mu_3= \max\{2\eta(2\mu_2), 6\mu_1\mu_2\}$,

\noindent $(3)$ $\mu_4=\max\{4^{16c'CM\mu_3(\eta(6\mu_1)+1)},
(\mu_3(\eta(6\mu_1)+1))^{16c'CM}\} $,

\noindent $(4)$ % % %$\mu_5=\max\Big\{\frac{16c'M\mu_3\mu_4}{\eta^{-1}\big(\frac{1}{4\mu_3\mu_4}\big)},16 %c'M\mu_3\mu_4\eta^{-1}\big(\frac{1}{4\mu_3\mu_4}\big)
% \Big\}$,
$\mu_5=  16c'M\mu_3\mu_4 \,\max \big\{1/u, u
 \big\}, \quad u = \eta^{-1}\big({1}/{(4\mu_3\mu_4)}\big) $,

\noindent $(5)$ $\mu_6=4\mu_4\mu_5$,

\noindent $(6)$ $\mu_7=\big(\mu_6\big(\eta(2\mu_2)+1\big)\big)^{8c'M}$,

\noindent where $\mu_1$ and $c'$ $(\leq 7c^3)$ are from Lemma
\ref{lem-j-j} and Theorem \Ref{thm0.1}, respectively, and $M$, $C$
and $\eta$ from {\it Basic assumption B}.

\subsection{Several lemmas}

\begin{lem} \label{lem-3-0'}  For $v_1\in \partial D$ and $v_2\in D$,
if $|v_1-v_2|\leq 2\mu
_2\,d_{D}(v_2)$, then
$$|v'_1-v'_2|\leq\mu_3\,d_{D'}(v'_2).$$
\end{lem}

\bpf
Let $x'_1\in \partial D'$ be such that $d_{D'}(v'_2)\geq
\frac{1}{2}|x'_1-v'_2|$.
It follows from the assumptions on $f$ that
$$\frac{|v'_1-v'_2|}{|x'_1-v'_2|}\leq \eta(2\mu_2),$$
whence $$|v'_1-v'_2|\leq 2\eta(2\mu_2)d_{D'}(v'_2),$$
which shows that the lemma holds.\epf

For $z'_1$, $z'_2\in D' \subset E'$, let $\gamma'$ be a
$2$-neargeodesic joining $z'_1$ and $z'_2$ in $D'$. In the following, we aim to prove that $\gamma'$ satisfies the conditions
\eqref{wx-4} and \eqref{wx-5} in Definition \ref{def1.3}.

Without loss of generality, we may assume that $d_{D'}(z'_2)\geq
d_{D'}(z'_1)$. Let $x_0\in\gamma$ be the first point in the
direction from $z_1$ to $z_2$ such that

$$d_{D}(x_0)=\sup\limits_{p\in \gamma}d_{D}(p).
$$

\begin{lem} \label{lem-3-1}  For all  $z'\in
\gamma'[z'_1, x'_0]$,
 $$\diam(\gamma'[z'_1, z'])\leq
\mu_4\;d_{D'}(z'),$$ and for all  $z'\in \gamma'[z'_2, x'_0]$,
$$\diam(\gamma'[z'_2, z'])\leq \mu_4\;d_{D'}(z').$$\end{lem}

\bpf We only need to prove the former assertion since the proof for
the latter one is similar. We prove it by a contradiction. Suppose
there exists some point $z'_{11}\in \gamma'[z'_1, x'_0]$ such that

\beq\label{w-1}\diam(\gamma'[z'_1, z'_{11}])>
\mu_4\;d_{D'}(z'_{11}).\eeq

Obviously, there exists some point $w'_{11}\in\gamma'[z'_1,
z'_{11}]$ such that $$|w'_{11}-z'_{11}|=
\frac{1}{2}\diam(\gamma'[z'_1, z'_{11}]).$$ Then, by (\ref{w-1}), we
see that $$|w'_{11}-z'_{11}|>\frac{\mu_4}{2}\;d_{D'}(z'_{11}),$$
whence by (\ref{eq(0000)})
\begin{eqnarray*}
 k_{D}(w_{11}, z_{11})&\geq& \frac{1}{M}\big(k_{D'}(w'_{11}, z'_{11})-C\big)\geq
  \frac{1}{M}\left (\log\Big(1+\frac{|w'_{11}-z'_{11}|}{d_{D'}(z'_{11})}\Big)-C\right )\\ \nonumber&\geq&
\frac{1}{M}\left(\log\Big(1+\frac{\mu_4}{2}\Big)
-C\right)
>1,\end{eqnarray*}
which,  together with (\ref{upperbdk}), implies that

\beq\label{l-h-w-1}|w_{11}-z_{11}|> \frac{1}{2}\max\{d_{D}(z_{11}),d_{D}(w_{11})\}.\eeq
Let $w_{12}\in \partial D$ be such that
$$|w_{12}-w_{11}|\leq 2d_{D}(w_{11}).$$
 Since $\mu_2\geq 1$, it follows from Lemma
\ref{lem-3-0'} that $$|w'_{12}-w'_{11}|\leq \mu_3 d_{D'}(w'_{11}).$$
Lemma \ref{lem-j-j} and (\ref{l-h-w-1}) imply
 \beq\label{w-3}|w_{12}-z_{11}|\leq |w_{12}-w_{11}|+|w_{11}-z_{11}|\leq
5\mu_1\,d_{D}(z_{11}).\eeq Therefore we see from Lemma \ref{lem-3-0'} and the
fact ``$5\mu_1\leq \mu_2$" that \beq\label{w-4} |w'_{12}-
z'_{11}|\leq \mu_3\,d_{D'}(z'_{11}).\eeq

On one hand, if $|w'_{12}-w'_{11}|\leq \mu_3\eta(6\mu_1)|w'_{12}-z'_{11}|$,
then, by (\ref{w-4}), we have
\begin{eqnarray*}
\diam(\gamma'[z'_1, z'_{11}])&=& 2|w'_{11}-z'_{11}|
\leq 2(|w'_{12}-w'_{11}|+|w'_{12}-z'_{11}|)\\ \nonumber &\leq&
2\big(\mu_3\eta(6\mu_1)+1\big)|w'_{12}- z'_{11}|\\
&\leq&  2\mu_3\big(\mu_3\eta(6\mu_1)+1\big)\,d_{D'}(z'_{11}),\end{eqnarray*} which  contradicts (\ref{w-1}).

On the other hand, if $|w'_{12}-w'_{11}|> \mu_3\eta(6\mu_1)|w'_{12}-z'_{11}|$, then
Lemma \ref{lem-j-j} and (\ref{w-3}) imply that
$$
|w_{12}-w_{11}|\leq |w_{12}-z_{11}|+|w_{11}-z_{11}|\leq 6\mu_1\,d_{D}(z_{11}).$$ Hence
$$\mu_3\eta(6\mu_1)\leq \frac{|w'_{12}-w'_{11}|}{|w'_{12}-z'_{11}|}\leq \eta(6\mu_1),$$
since $\frac{|w_{12}-w_{11}|}{|w_{12}-z_{11}|}\leq 6\mu_1$. This is the desired contradiction.\epf

\begin{lem} \label{lem-3-2}  For all  $z'\in
\gamma'[z'_1, z'_2]$, we have $\ds\min_{j=1,2}\ell (\gamma' [z'_j, z'])\leq
\mu_6\,d_{D'}(z')$.\end{lem}

\bpf We use $z'_0\in\gamma'$ to denote the
first point on $\gamma'$ in the direction from $z'_1$ to $z'_2$ such that

\be\label{wx-1} d_{D'}(z'_0)=\sup\limits_{p'\in \gamma'}d_{D'}(p').
\ee
It is possible that $z'_0=z'_1$ or $z'_2$. Obviously, there exists a
nonnegative integer $m$ such that

$$ 2^{m}\, d_{D'}(z'_1) \leq d_{D'}(z'_0)< 2^{m+1}\, d_{D'}(z'_1),
$$ and we use $y'_0$ to denote the first point on $\gamma[z'_1,z'_0]$
from $z'_1$ to $z'_0$ with

$$d_{D'}(y'_0)=2^{m}\, d_{D'}(z'_1).
$$

We define a sequence $\{x'_k\}$. Let $x_1'=z_1'$. If
$y'_0=z'_1\not=z'_0$, then we let $x'_2=z'_0$. If $y'_0\not= z'_1$,
then we let $x'_2$, $\ldots$, $x'_{m+1}\in \gamma'[z'_1,z'_0]$ be
the points such that for each $i\in \{2$, $\ldots,m+1\}$, $x'_i$
is the first point from $z'_1$ to $z'_0$ with
$$d_{D'}(x'_i)=2^{i-1}\, d_{D'}(x'_1).
$$
Therefore $x'_{m+1}=y'_0$. In the case $y'_0\not= z'_0$,  we define $x'_{m+2}=z'_0$.

 Similarly, let $s\geq 0$ be the integer such that

$$ 2^{s}\, d_{D'}(z'_2) \leq d_{D'}(z'_0)< 2^{s+1}\, d_{D'}(z'_2),
$$
and let $x_{1,0}'$ be the first point on $\gamma'[z'_2,z'_0]$ from
$z'_2$ to $z'_0$ with

$$d_{D'}(x_{1,0}')=2^{s}\, d_{D'}(z'_2).
$$

In the same way, we define another sequence $\{x'_{1,k}\}$. We let
$x_{1,1}'=z_2'$. If $x_{1,0}'= z'_2\not=z'_0$, then we let
$x_{1,2}'=z'_0$. If $x_{1,0}'\not= z'_2$, then we let $x_{1,2}'$,
$\ldots$, $x_{1,s+1}'$ be the points on $\gamma'[z'_2,z'_0]$ such
that $x_{1,j}'$ $(j=2$, $\ldots$, $s+1)$ denotes the first point
from $x_{1,1}'$ to $z'_0$ with
%\be\label{hws-eq(4.3)}
$$d_{D'}(x_{1,j}')=2^{j-1}\, d_{D'}(x_{1,1}').
$$
Then $x_{1,s+1}'=x'_{1,0}$. If $x_{1,0}'\not= z'_0$, then we let
$x_{1,s+2}'=z'_0$.
\medskip

For a proof of the lemma, it is enough to prove that for every $z'\in\gamma'[z'_1,z'_0]$,

\be\label{wx-2} \ell(\gamma'[z'_1,z'])\leq \mu_6\,d_{D'}(z'),\ee and for all $z'\in\gamma'[z'_2,z'_0]$,

\be\label{wx-2'} \ell(\gamma'[z'_2,z'])\leq \mu_6\,d_{D'}(z').\ee

We only need to prove \eqref{wx-2} since the proof for \eqref{wx-2'} is similar.

Before the proof of \eqref{wx-2}, we prove two claims.

\bcl\label{cl--1-2}For each $i\in \{1,\cdots,m+1\}$, if $y'\in
\gamma'[x'_i,x'_{i+1}]$, then $d_{D'}(x'_i)\leq 2\mu_4\,d_{D'}(y')$.
\ecl

To prove this claim, it suffices to consider the case: $d_{D'}(y')<
\frac{1}{2}d_{D'}(x'_i)$ since the proof for the case:
$d_{D'}(y')\geq \frac{1}{2}d_{D'}(x'_i)$ is trivial $(\mu_4\geq 1)$.
It follows from $|x'_i-y'|\geq d_{D'}(x'_i)-d_{D'}(y')$ that

$$\min\{|x'_i-y'|, |y'-x'_{i+1}|\}>\frac{1}{2}d_{D'}(x'_i).$$

If $\gamma'[z'_1,y']\subset \gamma'[z'_1,x'_0]$, then, by Lemma
\ref{lem-3-1}, $$d_{D'}(y')\geq \frac{1}{\mu_4}|x'_i-y'|\geq
\frac{1}{2\mu_4}d_{D'}(x'_i).$$
For the remaining case, we know that $\gamma'[z'_2,y']\subset \gamma'[z'_2,x'_0]$. The similar
reasoning as above shows that $$d_{D'}(y')\geq
\frac{1}{\mu_4}|y'-x'_{i+1}|\geq \frac{1}{2\mu_4}d_{D'}(x'_i).$$
Hence the proof of Claim \ref{cl--1-2} is complete.

\bcl\label{cl--1}For all $i\in \{1,\cdots,m+1\}$,
$\ell(\gamma'[x'_i,x'_{i+1}])\leq \mu_5\,d_{D'}(x'_i)$. \ecl

Suppose on the contrary that there exists some $i\in
\{1,\cdots,m+1\}$ such that
\be\label{h-4}\ell(\gamma'[x'_i,x'_{i+1}])> \mu_5\,d_{D'}(x'_i).\ee
Because $\gamma'$ is a $2$-neargeodesic, we get by (\ref{h-4})
$$k_{D'}(x'_i,x'_{i+1})\geq
\frac{1}{2}\int_{\gamma'[x'_i,x'_{i+1}]}\frac{|dx'|}{d_{D'}(x')}>\frac{\mu_5}{4}\,,$$
and we see that

$$ k_{D}(x_i,x_{i+1})\geq
\frac{1}{M}k_{D'}(x'_i,x'_{i+1})-\frac{C}{M}> 2c'\,\mu_3\,\mu_4 \,
\max\left\{1/u,u  \right\},$$ where $ u =
\eta^{-1}\big({1}/{(4\mu_3\mu_4)}\big)$. Then it follows from the
inequality:
$$k_{D}(x_i,x_{i+1})\leq
c'\log\left(1+\frac{|x_i-x_{i+1}|}{\min\{d_{D}(x_i),d_{D}(x_{i+1})\}}\right)$$ that

\be\label{h-6'}\min\{d_{D}(x_i),d_{D}(x_{i+1})\}\leq
\min\left\{\frac{1}{e^{ 2\mu_3\mu_4/ u }-1},
\frac{1}{e^{2\mu_3\mu_4}-1}\right\}|x_i-x_{i+1}|.\ee
Without loss of generality, we may assume that

\be\label{w'-1}\min\{d_{D}(x_i),d_{D}(x_{i+1})\}=d_{D}(x_i).\ee
Take $x_{2i}\in \partial D$ such that

\be\label{w'-2} |x_{2i}-x_i|\leq 2d_{D}(x_i).\ee
Then Lemma \ref{lem-3-0'} implies

\be\label{w'-3}|x'_{2i}-x'_i|\leq
\mu_3d_{D'}(x'_i).\ee

If $\gamma'[z'_1,x'_{i+1}]\subset \gamma'[z'_1,x'_0]$ or
$\gamma'[z'_2,x'_i]\subset \gamma'[z'_2,x'_0]$,  then, by Lemma
\ref{lem-3-1}, $$|x'_i-x'_{i+1}|\leq \mu_4d_{D'}(x'_{i+1}) \leq 2\mu_4d_{D'}(x'_i).$$
For the remaining case, we know that $x'_0\in\gamma'[x'_i,x'_{i+1}]$. Then Lemma \ref{lem-3-1} yields
$$|x'_i-x'_{i+1}| \leq 2\max\{|x'_i-x'_0|, |x'_{i+1}-x'_0|\} \leq 2\mu_4d_{D'}(x'_0)  \leq 4\mu_4d_{D'}(x'_i).$$
Hence we get that for each $i\in \{1,\cdots, m+1\}$,

\be\label{l-10}|x'_i-x'_{i+1}|\leq 4\mu_4d_{D'}(x'_i).\ee
By (\ref{w'-3}) and (\ref{l-10}), we have

\be\label{l-10'}
|x'_{2i}-x'_{i+1}| \leq
|x'_i-x'_{i+1}|+|x'_{2i}-x'_i|
\leq (\mu_3+4\mu_4)d_{D'}(x'_i),\ee and by (\ref{h-6'}), (\ref{w'-1}) and (\ref{w'-2}),
$$
|x_{2i}-x_{i+1}| \geq
|x_i-x_{i+1}|-|x_{2i}-x_i|\geq
\max \left\{ \frac{e^{ 2\mu_3\mu_4/u }}{4},
\frac{e^{2\mu_3\mu_4}-3}{2} \right\}|x_{2i}-x_i|.$$
Hence \eqref{l-10'} implies

$$ \frac{1}{\mu_3+4\mu_4}\leq
\frac{|x'_i-x'_{2i}|}{|x'_{2i}-x'_{i+1}|} \leq
\eta\left( u \right)
= \frac{1}{4\mu_3\mu_4}.$$ This is the desired
contradiction, which completes the proof of Claim \ref{cl--1}.
\medskip

Now we are ready to prove \eqref{wx-2}.\medskip

If $z'\in \gamma'[z'_1, x'_{m+1}]$, then there exists some $k\in
\{1,\cdots,m\}$ such that $z'\in \gamma'[x'_k, x'_{k+1}]$. If $k=1$,
then it easily follows from Claims \ref{cl--1-2} and \ref{cl--1}
that

\beq\label{l-11} \ell(\gamma'[z'_1,z'])\leq
\ell(\gamma'[x'_1,x'_2])\leq \mu_5d_{D'}(x'_1)\leq
2\mu_4\mu_5d_{D'}(z').\eeq  If $k>1$, then, again, by Claims  \ref{cl--1-2} and \ref{cl--1},

\beq\label{l-12}
\ell(\gamma'[z'_1,z'])&\leq& \ell(\gamma'[x'_1,x'_2])+\cdots + \ell(\gamma'[x'_{k-1},x'_k])+\ell(\gamma'[x'_k,z'])\\
\nonumber &\leq& \mu_5\big(d_{D'}(x'_1)+\cdots+d_{D'}(x'_{k-1})+d_{D'}(x'_k)\big)\\
\nonumber &\leq&2\mu_5d_{D'}(x'_k)\\
\nonumber &\leq& 4\mu_4\mu_5d_{D'}(z').\eeq

Now we consider the remaining case: $z'\in \gamma'[x'_{m+1}, z'_0 ]$. We
infer from Claims   \ref{cl--1-2} and \ref{cl--1} that

\beq\label{l-13}
\ell(\gamma'[z'_1,z'])&\leq& \mu_5\big(d_{D'}(x'_1)+d_{D'}(x'_2)+\cdots+d_{D'}(x'_{m})+d_{D'}(x'_{m+1})\big) \\
\nonumber &\leq& 2\mu_5d_{D'}(x'_{m+1})\\ \nonumber &\leq& 4\mu_4\mu_5d_{D'}(z').\eeq
The combination of \eqref{l-11}, \eqref{l-12} and \eqref{l-13} shows that
 for all
$z'\in \gamma'[z'_1, z'_0]$, $$\ell(\gamma'[z'_1,z'])\leq
4\mu_4\mu_5d_{D'}(z').$$
Hence  Lemma \ref{lem-3-2} holds.\epf

Further, we have

\begin{lem} \label{lem-3-3} $\ell (\gamma' [z'_1, z'_2])\leq
\mu_7\,|z'_1-z'_2|$.\end{lem}

\bpf Suppose, on the contrary, that
\be\label{h-7}\ell (\gamma' [z'_1, z'_2])>
\mu_7\,|z'_1-z'_2|.\ee
We first prove a claim.

\bcl\label{h-8} $d_{D'}(z'_2)\leq 7\,|z'_1-z'_2|$. \ecl

Also we prove this claim by contradiction.
Suppose

\be\label{hmw} d_{D'}(z'_2)> 7\,|z'_1-z'_2|.\ee
Because $\gamma'$ is a $2$-neargeodesic, we have by
(\ref{eq(0000)}) that

\begin{eqnarray*}\frac{1}{2}\log\Big(1+\frac{\ell (\gamma' [z'_1, z'_2])}{d_{D'}(z'_1)}\Big)
&\leq&\frac{1}{2}\ell_k(\gamma' [z'_1, z'_2])\leq
k_{D'}(z'_1,z'_2)\leq \int_{[z'_1,z'_2]}\frac{|dz'|}{d_{D'}(z')}\\ \nonumber   &\leq&
\frac{7}{6}\cdot\frac{|z'_1-z'_2|}{d_{D'}(z'_2)}< \frac{1}{6},\end{eqnarray*}
since $d_{D'}(z')\geq d_{D'}(z'_2)-|z'_2-z'_1|$ for all $z'\in [z'_1,z'_2]$.
By \eqref{h-7} and \eqref{hmw}, this is a contradiction. Hence Claim \ref{h-8} holds true.
\medskip

Recall that $z'_0\in \gamma'$ satisfies \eqref{wx-1}. Let
$x'$ be the point of $\gamma'$ which bisects the arclength of $\gamma'\,,$ i.e.
 $ \ell(\gamma'[z_1', x']) = \ell(\gamma'[z_2', x']) \,.$ Then Lemma \ref{lem-3-2} implies
$$ \ell(\gamma'[z_1', x'])\leq \mu_6d_{D'}(x')\leq
\mu_6d_{D'}(z_0').$$ Hence it follows from (\ref{h-7}) and Claim
\ref{h-8} that

\be\label{h-9} d_{D'}(z'_2)\leq 7|z'_1-z'_2|< \frac{7}{\mu_7}\ell
(\gamma' [z'_1, z'_2])\leq  \frac{14\mu_6}{\mu_7}d_{D'}(z'_0),\ee
whence by Theorem \Ref{thm0.1}
\begin{eqnarray*}c'\log\Big(1+\frac{|z_2-z_0|}{\min\{d_{D}(z_2),d_{D}(z_0)\}}\Big)&\geq&
k_{D}(z_2,z_0)\geq \frac{1}{M}k_{D'}(z'_2,z'_0)-\frac{C}{M}\\
&\geq&
\frac{1}{M}\log\frac{d_{D'}(z'_0)}{d_{D'}(z'_2)}-\frac{C}{M}\\ &>&
\frac{1}{2M}\log \frac{\mu_7}{\mu_6}\end{eqnarray*} and
\begin{eqnarray*}c'\log\Big(1+\frac{|z_1-z_0|}{\min\{d_{D}(z_1),d_{D}(z_0)\}}\Big)&\geq&
k_{D}(z_1,z_0)\geq \frac{1}{M}k_{D'}(z'_1,z'_0)-\frac{C}{M}\\
&\geq&
\frac{1}{M}\log\frac{d_{D'}(z'_0)}{d_{D'}(z'_1)}-\frac{C}{M}\\&>&
\frac{1}{2M}\log \frac{\mu_7}{\mu_6}.\end{eqnarray*} These show that

\be\label{h-10}\min\{d_{D}(z_2),d_{D}(z_0)\}<
\frac{1}{\big(\mu_6\big(\eta(2\mu_2)+1\big)\big)^2}\,|z_0-z_2|,\ee
$$\min\{d_{D}(z_1),d_{D}(z_0)\}<
\frac{1}{\big(\mu_6\big(\eta(2\mu_2)+1\big)\big)^2}\,|z_0-z_1|,$$
and we see by (\ref{upperbdk}) \be\label{h-6} |z'_2- z'_0|>
\frac{1}{2}d_{D'}(z'_0) \;\;\mbox{and}\;\; |z'_1- z'_0|>
\frac{1}{2}d_{D'}(z'_0), \ee since $k_{D'}(z_1',z_0')> 1$ and
$k_{D'}(z_2',z_0')> 1$.

\bcl\label{h-15} $\min\{d_{D}(z_2),d_{D}(z_1)\}\leq 2\mu_1d_{D}(z_0)$.\ecl
\medskip

We first prove that
$d_D(z_1)\leq 2\mu_1d_D(z_0)$ when $\gamma[z_1,z_0]\subset\gamma[z_1,x_0]$.

If $|z_1-z_0|\leq \frac{1}{2}d_D(z_1)$, then $$d_D(z_0)\geq d_D(z_1)-|z_1-z_0|\geq \frac{1}{2}d_D(z_1).$$

On the other hand, if $|z_1-z_0|> \frac{1}{2}d_D(z_1)$, then we obtain from Lemma \ref{lem-j-j} that
$$d_D(z_1)\leq 2|z_1-z_0|\leq 2\mu_1d_D(z_0).$$

A similar discussion as above shows that $d_D(z_2)\leq 2\mu_1d_D(z_0)$ when $\gamma[z_2,z_0]\subset\gamma[z_2, x_0]$.
The proof of Claim \ref{h-15} is complete.
\medskip

Without loss of generality, we may assume that
$$\min\{d_{D}(z_2),d_{D}(z_1)\}=d_{D}(z_2).$$ Then (\ref{h-10}) and Claim \ref{h-15} imply

\be\label{h-10'}d_{D}(z_2)\leq
\frac{1}{\mu_6\big(\eta(2\mu_2)+1\big)}\,|z_0-z_2|.\ee
Take $w_{13}\in\partial D$ such that

\be\label{wx-3} |w_{13}-z_2|\leq 2d_{D}(z_2).\ee  It follows from Lemma \ref{lem-3-0'} and Claim \ref{h-8}  that
 \be\label{h-11'}|w'_{13}-z'_2|\leq \mu_3d_{D'}(z'_2)\leq
7\mu_3|z'_1-z'_2|,\ee whence \be\label{h-12}|w'_{13}-z'_1|\leq
|w'_{13}-z'_2|+|z'_2-z'_1|\leq (1+7\mu_3)|z'_1-z'_2|.\ee

By the $2$-neargeodesic property of $\gamma'\,,$ (\ref{eq(0000)}), \eqref{h-7} and Claim \ref{h-8}, we have

\begin{eqnarray*}k_D(z_1,z_2)&\geq&
 \frac{1}{M}k_{D'}(z'_1,z'_2)-\frac{C}{M}\geq \frac{1}{2M}\ell_k(\gamma' [z'_1,
z'_2])-\frac{C}{M}\\
&\geq&\frac{1}{2M}\log\Big(1+\frac{\ell (\gamma' [z'_1,
z'_2])}{d_{D'}(z'_1)}\Big)-\frac{C}{M}\geq \frac{1}{4M}\log
\mu_7 \\&>&1,\end{eqnarray*} which shows by (\ref{upperbdk})
\be\label{h-10''}|z_1-z_2|\geq \frac{1}{2}d_{D}(z_1).\ee

Because of \eqref{h-10'} and \eqref{wx-3}, we get
$$|w_{13}-z_0|\leq |w_{13}-z_2|+|z_2-z_0|\leq
\Big(1+\frac{2}{\mu_6\big(\eta(2\mu_2)+1\big)}\Big)|z_2-z_0|.$$
Next by Lemma \ref{lem-j-j} and \eqref{h-10''}, we have

\beq\label{Tu-2}
|z_0-z_2|&\leq &\diam(\gamma) \leq \mu_1\max\Big\{|z_1-z_2|, 2\Big(e^{h}-1\Big)d_{D}(z_2)\Big\}\\ \nonumber
&\leq & \mu_1\max\Big\{1, 4\Big(e^{h}-1\Big)\Big\}|z_1-z_2|\\ \nonumber &\leq & \mu_2|z_1-z_2|.
\eeq
Hence (\ref{h-10'}), \eqref{wx-3} and (\ref{Tu-2})  yield

\beq\label{h-11}|w_{13}-z_1|&\geq& |z_1-z_2|-|w_{13}-z_2|\\
\nonumber&\geq&
\Big(\frac{1}{\mu_2}-\frac{2}{\mu_6\big(\eta(2\mu_2)+1\big)}\Big)\,|z_0-z_2|\\
\nonumber &\geq& \frac{1}{2\mu_2}|w_{13}-z_0|.\eeq We see from
(\ref{h-9}), (\ref{h-6}) and (\ref{h-11'}) that

\beq\label{h-13}|w'_{13}-z'_0|&\geq& |z'_2-z'_0|-|w'_{13}-z'_2|\geq
\Big(\frac{\mu_7}{4\mu_6}-7\mu_3\Big)\,|z'_1-z'_2|\\
\nonumber &>&  \mu_6\big(\eta(2\mu_2)+1\big)|z'_1-z'_2|.\eeq The combination of
 (\ref{h-12}), (\ref{h-11}) and (\ref{h-13}) shows
$$\eta(2\mu_2)<
\frac{|w'_{13}-z'_0|}{|w'_{13}-z'_1|}
\leq  \eta(2\mu_2).$$ This desired
contradiction completes the proof of Lemma \ref{lem-3-3}.
\epf

The next result is obvious from a similar argument as above.

\bcor\label{fri-1} Under the assumptions of Theorem \ref{thm1.1}, if
$f$ extends to a homeomorphism $\overline{f}: \overline{D}\to
\overline{D}'$ and $\overline{f}$ is $\eta$-QM rel $\partial D$,
then each $a$-neargeodesic $(a>1)$ in $D'$ is a double $a'$-cone
arc, where $a'$ depends only on $a$, $c$, $M$, $C$ and $\eta$. In
particular, if $a=2$, we can take $a'=\mu_7$. \ecor

\subsection{ The proof of
Theorem \ref{thm1.1}.}

For $z'_1$, $z'_2\in D'$, let $\gamma'$ be a $2$-neargeodesic
joining $z'_1$ and $z'_2$ in $D'$. Corollary \ref{fri-1}
 shows that $\gamma'$ is a double $\mu_7$-cone arc. Hence $D'$ is a $\mu_7$-uniform domain, and so the proof of
Theorem \ref{thm1.1} easily follows from Theorem \Ref{ThmC}. \qed


\begin{thebibliography}{99}

%\bibitem{Ahl}  {\sc L. V. Ahlfors}, Quasiconformal reflections, \textit{Acta Math.,} {\bf 109} (1963),
%291--301.

%\bibitem{FW}  {\sc F. W. Gehring}, Uniform domains and the ubiquitous
%quasidisks, \textit{Jahresber. Deutsch. Math. Verein,} {\bf
%89}(1987), 88--103.

\bibitem{Avv}  {\sc G. D. Anderson, M. K. Vamanamurthy and M. Vuorinen}, Dimension-free
quasiconformal distortion in $n$-space. \textit{Trans. Amer. Math.
Soc.,} {\bf 297} (1986), 687--706.

\bibitem{Bea} {\sc A. F. Beardon},
The Apollonian metric of a domain in $\IR^{n}$,  In: Quasiconformal
mappings and analysis, Springer-Verlag, 1998, 91--108.

\bibitem{BHX}  {\sc
St. M. Buckley, D. A. Herron and X. Xie}, Metric space inversions,
quasihyperbolic distance, and uniform spaces, \textit{ Indiana Univ.
Math. J.,} {\bf 57} (2008), 837--890.


\bibitem{Geo}  {\sc F. W. Gehring and B. G. Osgood}, Uniform domains
 and the
quasi-hyperbolic metric, \textit{J. Analyse Math.,} {\bf 36} (1979),
50--74.

\bibitem{GP}  {\sc F. W. Gehring and B. P. Palka}, Quasiconformally
homogeneous domains, \textit{J. Analyse Math.,} {\bf 30} (1976),
172--199.

%\bibitem %[GO]
%{Ger6}  {\sc F. W. Gehring}, Injectivity of local quasi-isometries,
% \textit{Comment. Math. Helv.,} {\bf 57} (1982),
%202--220.


%\bibitem %[GO]
%{Ger5}  {\sc F. W. Gehring}, Extension of quasiisometric embeddings
%of Jordan curves, \textit{Complex Variables}, {\bf 5} (1986),
%245--263.

\bibitem {HIMPS} {\sc P.  H\"ast\"o, Z. Ibragimov, D. Minda, S. Ponnusamy and S. Sahoo},
{Isometries of some hyperbolic-type path metrics, and the hyperbolic medial axis.}  In the tradition of Ahlfors-Bers. IV,
 63--74, Contemp. Math., 432, Amer. Math. Soc., Providence, RI,  2007.


%\bibitem{Ml}  {\sc M. Huang and Y. Li}, Decomposition  properties of John domains in real normed vector
%spaces,  \textit{J. Math. Anal. Appl.,} {\bf 388}  (2012), 191--197.


\bibitem{k}  {\sc R. Kl\'en,} Local convexity properties of quasihyperbolic balls
in punctured space, \textit{J. Math. Anal. Appl.,} {\bf 342} (2008),
192--201.

\bibitem{krt}  {\sc R. Kl\'en, A. Rasila, and J. Talponen,}
Quasihyperbolic geometry in euclidean and Banach spaces, Proc.
ICM2010 Satellite Conf. International Workshop on Harmonic
 and Quasiconformal Mappings (HMQ2010), eds. D. Minda, S. Ponnusamy, N. Shanmugalingam,
 \textit{ J. Analysis,} {\bf 18} (2010), 261--278,
arXiv: 1104.3745v1 [math.CV].




%\bibitem{Yli}  {\sc Y. Li and X. Wang}, Unions of John domains and uniform domains in real normed vector
%spaces, \textit{Ann. Acad. Sci. Fenn. Ser. A I Math.}, {\bf 35}
%(2010), 627--632.

\bibitem{KN}  {\sc P. Koskela and T. Nieminen}, Uniform continuity of quasiconformal mappings and
conformal deformations, \textit{ Conform. Geom. Dyn.,} {\bf 12}
(2008), 10--17.

\bibitem{Martio-80}  {\sc O. Martio},
Definitions of uniform domains, \textit{Ann. Acad. Sci. Fenn. Ser. A
I Math.,} {\bf 5} (1980), 197--205.

\bibitem{rt} {\sc A. Rasila and J. Talponen,} Convexity properties of quasihyperbolic balls on
Banach spaces, \textit{Ann. Acad. Sci. Fenn.} 37, 2012, 215-228, arXiv:
1007.3197v1 [math. CV].

%\bibitem{Sc}  {\sc J. J. Sch\"affer}, Inner diameter, perimeter, and
%girth of spheres. \textit{Math. Ann.,} {\bf 173}(1967), 59-82.

\bibitem{TV}  {\sc P. Tukia and J. V\"{a}is\"{a}l\"{a}}, Quasisymmetric embeddings of metric spaces,
\textit{Ann. Acad. Sci. Fenn. Ser. A I Math.,} {\bf 5} (1980),
97-114.

%\bibitem{TV1}  {\sc P. Tukia and J. V\"{a}is\"{a}l\"{a}}, Bilipschitz extensions of maps having quasiconformal
%extensions, \textit{Math. Ann.,} {\bf 269} (1984), 651-572.

%\bibitem{Vai1}  {\sc J. V\"{a}is\"{a}l\"{a}}, Lectures on n-dimensional
%quasiconformal mappings, Springer-Verlag, 1971.

\bibitem{Vai2}  {\sc J. V\"{a}is\"{a}l\"{a}}, Quasim\"{o}bius maps,
\textit{J. Analyse Math.,} {\bf 44} (1985), 218--234.

\bibitem{Vai}  {\sc J. V\"{a}is\"{a}l\"{a}}, Uniform domains,
\textit{Tohoku Math. J.,} {\bf 40} (1988), 101--118.

\bibitem{Vai6-0}  {\sc J. V\"{a}is\"{a}l\"{a}}, Free quasiconformality
in Banach spaces. I, \textit{Ann. Acad. Sci. Fenn. Ser. A I Math.,}
{\bf 15} (1990), 355-379.

\bibitem{Vai6}  {\sc J. V\"{a}is\"{a}l\"{a}}, Free quasiconformality
in Banach spaces. II, \textit{Ann. Acad. Sci. Fenn. Ser. A I Math.,}
{\bf 16} (1991), 255-310.

\bibitem{Vai4}  {\sc J. V\"{a}is\"{a}l\"{a}}, Relatively and inner
 uniform domains,
\textit{Conformal Geom. Dyn.,} {\bf 2} (1998), 56--88.

%\bibitem{Vai7}  {\sc J. V\"{a}is\"{a}l\"{a}}, Free quasiconformality
%in Banach spaces. IV, \textit{Analysis and Topology, 697-717, World
%Sci. Publ., River Edge, N. J.,} 1998.

\bibitem{Vai5} {\sc J. V\"{a}is\"{a}l\"{a}}, The free quasiworld,
freely quasiconformal and related maps in Banach spaces,
\textit{Banach Center Publ.,} {\bf 48} (1999), 55--118.

\bibitem{Vai3}  {\sc J. V\"{a}is\"{a}l\"{a}}, Quasihyperbolic
geodesics in convex domains,
\textit{Results Math.,} {\bf 48} (2005), 184--195.



\bibitem{Vai9} {\sc J.
 V\"{a}is\"{a}l\"{a},} Quasihyperbolic geometry of planar domains. \textit{Ann. Acad. Sci. Fenn. Math.} {\bf 34} (2009),  447--473.




\bibitem{Mvo1} {\sc M. Vuorinen}, Conformal invariants and
quasiregular mappings, \textit{J. Analyse Math.,} {\bf 45} (1985),
69--115.

\end{thebibliography}
\end{document}